\documentstyle[10pt]{article}
\begin{document}
\begin{center}
{\Large \bf AF-domains and their generalizations }\vspace{1mm}
\end{center}
\begin{center}
{\bf Samir Bouchiba}\\
{\footnotesize {\it Department of Mathematics, University Moulay
Ismail, Meknes 50000, Morocco}} \vspace{8mm}
\end{center}

\noindent {\large \bf
------------------------------------------------------------------------------------------}
\noindent {\footnotesize \bf Abstract}\\

{\footnotesize In this paper, we are concerned with the study of the
dimension theory of tensor products of algebras over a field $k$. We
introduce and investigate the notion of generalized AF-domain
(GAF-domain for short) and prove that any $k$-algebra $A$ such that
the polynomial ring in one variable $A[X]$ is an AF-domain is in
fact a GAF-domain, in particular any AF-domain is a GAF-domain.
Moreover, we compute the Krull dimension of $A\otimes_kB$ for any
$k$-algebra $A$ such that $A[X]$ is an AF-domain and any $k$-algebra
$B$ generalizing the main
theorem of Wadsworth in [16].}\\

\noindent {\footnotesize {\it MSC} (2000): 13C15; 13B24.}\\
\noindent {\large \bf
------------------------------------------------------------------------------------------}
\vspace{5mm}

\noindent {\bf 1. Introduction}\\

All rings considered in this paper are commutative with identity
element and all ring homomorphisms are unital. Throughout, $k$
stands for a field. We shall use t.d.($A:k)$, or t.d.($A)$ when no
confusion is likely, to denote the transcendence degree of a
$k$-algebra $A$ over $k$ (for nondomains, t.d.$(A):=$ sup$\Big
{\{}$t.d.$(\displaystyle {\frac Ap}):p\in$ Spec$(A)\Big {\}}$),
$A[n]$ to denote the polynomial ring $A[X_1,...,X_n]$ and $p[n]$ to
denote the prime ideal $p[X_1,...,X_n]$ of $A[X_1,...,X_n]$ for each
prime ideal $p$ of $A$. Also, we use Spec($A)$ to denote the set of
prime ideals of a ring $A$ and $\subset$ to denote proper set
inclusion. All $k$-algebras considered throughout this paper are
assumed to be of finite transcendence degree over $k$. Any
unreferenced material is standard as in [8], [12], [13] and [14].

Several authors have been interested in studying the prime ideal
structure and related topics of tensor products of algebras over a
field $k$. The initial impetus for these investigations was a paper
of R. Sharp on Krull dimension of tensor products of two extension
fields. In fact, in [15], Sharp proved that, for any two extension
fields $K$ and $L$ of $k$, dim$(K\otimes_kL)=$ min(t.d.($K),$
t.d.($L))$ (actually, this result appeared ten years earlier in
Grothendieck's EGA [10, Remarque 4.2.1.4, p. 349]). This formula is
rather surprising since, as one may expect, the structure of the
tensor product should reflect the way the two components interact
and \\
\noindent {\large \bf ---------------}

{\footnotesize {\it E-mail address}: sbouchiba@hotmail.com

{\it keywords:} Krull dimension, tensor product, prime ideal,
AF-domain.}\newpage

\noindent not only the structure of each component. This fact is
what most motivated Wadsworth's work in [16] on this subject. His
aim was to seek geometrical properties of primes of $A\otimes_kB$
and to widen the scope of algebras $A$ and $B$ for which
dim$(A\otimes_kB)$ depends only on individual characteristics of $A$
and $B$. The algebras which proved tractable for Krull dimension
computations turned out to be those domains $A$ which satisfy the
altitude formula over $k$ (AF-domains for short), that is,
$$ht(p)+\mbox {t.d.}(\frac Ap)=\mbox { t.d.}(A)$$ for all prime
ideals $p$ of $A$. It is worth noting that the class of AF-domains
contains the most basic rings of algebraic geometry, including
finitely generated $k$-algebras that are domains. Wadsworth proved,
via [16, Theorem 3.8], that if $A_1$ and $A_2$ are AF-domains, then
$$\mbox {dim}(A_1\otimes_kA_2)=\mbox { min}\Big {(}\mbox
{dim}(A_1)+\mbox {t.d.}(A_2),\mbox { t.d.}(A_1)+\mbox {dim}(A_2)\Big
{)}.$$ \noindent His main theorem stated a formula for
dim($A\otimes_kB)$ which holds for an AF-domain $A$, with no
restriction on $B$, namely:$$ \begin{array}{lll} \mbox
{dim}(A\otimes_kB)&=&D\Big {(}\mbox {t.d.}(A),\mbox {dim}(A),B\Big
{)}\\
&:=&\mbox {max}\Big {\{}ht(q[\mbox {t.d.}(A)])+\mbox {min}\Big
{(}\mbox {t.d.}(A),\mbox { dim}(A)+\displaystyle {\mbox {t.d.}(\frac
Bq})\Big {)}:q\in \mbox {Spec}(B)\Big {\}}\\
&&\mbox {[16, Theorem 3.7].}
\end{array}$$

On the other hand, in [11], Jaffard proved that, for any ring $A$
and any positive integer $n$, the Krull dimension of $A[n]$ can be
realized as the length of a special chain of $A[n]$. Recall that a
chain $C=\{Q_0\subset Q_1\subset ...\subset Q_s\}$ of prime ideals
of $A[n]$ is called a special chain if for each $Q_i$, the ideal
$(Q_i\cap A)[n]$ belongs to $C$. Subsequently, based on the thorough
and brilliant work of J. Arnold in [1], Brewer et al. gave an
equivalent and simple version of Jaffard's theorem. Actually, they
showed that, for each positive integer $n$ and each prime ideal $P$
of $A[n]$, $ht(P)=ht(q[n])+ht(\displaystyle {\frac P{q[n]})}$ [7,
Theorem 1], where $q:=P\cap A$. Taking into account the natural
isomorphism $B[n]\cong k[n]\otimes_kB$ for each $k$-algebra $B$, we
generalized in [5] this special chain theorem to tensor products of
$k$-algebras. Effectively, we proved that if $A$ and $B$ are
$k$-algebras such that $A$ is an AF-domain, then for each prime
ideal $P$ of $A\otimes_kB$,
$$ht(P)=ht(A\otimes_kq)+ht(\frac P{A\otimes_kq})=ht(q[\mbox
{t.d.}(A)])+ht(\frac P{A\otimes_kq}),$$ where $q=P\cap B$  (cf. [5,
Lemma 1.5]). It turns out that this very geometrical property
totally characterizes the AF-domains. In fact, we proved, in [4],
that the following statements are equivalent for a domain $A$ which
is a $k$-algebra:\newpage

a) $A$ is an AF-domain;

b) $A$ satisfies SCT (for special chain theorem), that is, for each
$k$-algebra $B$ and each prime ideal $P$ of $A\otimes_kB$ with
$q:=P\cap B$,
$$ht(P)=ht(q[\mbox {t.d.}(A)])+ht(\displaystyle {\frac
P{A\otimes_kq})}=ht(A\otimes_kq)+ht(\frac P{A\otimes_kq}) \mbox {
[4, Theorem 1.1]}.$$ In view of this, it is then natural to
generalize the AF-domain notion by setting the following
definitions:

We say that a $k$-algebra $A$ satisfies GSCT (for generalized
special chain theorem) with respect to a $k$-algebra $B$ if
\begin{center}
$ht(P)=ht(p\otimes_kB+A\otimes_kq)+ht(\displaystyle {\frac
P{p\otimes_kB+A\otimes_kq})}$\\

 for each prime ideal $P$ of
$A\otimes_kB$, with $p=P\cap A$ and $q=P\cap B$,
\end{center}

\noindent and we call a generalized AF-domain (GAF-domain for short)
a domain $A$ such that $A$ satisfies GSCT with respect to any
$k$-algebra $B$.\\ There is no known example in the literature of a
$k$-algebra $A$ which is a domain and which is not a GAF-domain.
This may lead one to ask whether any $k$-algebra which is a domain
is a GAF-domain. The
object of this paper is to handle the following question:\\

(Q): Is any domain $A$ which is a $k$-algebra such that the
polynomial ring $A[n]$ is an AF-domain, for some positive integer
$n$, a GAF-domain?\\

\noindent It is significant, in this regard, to note that if $A$ is
an AF-domain then $A[n]$ is an AF-domain for each integer $n\geq 0$,
and using pullback constructions, Proposition 2.1 shows that these
two notions do not coincide by providing a family of $k$-algebras
$A$ such that $A$ is not an AF-domain while there exists a positive
integer $r$ such that the polynomial ring $A[r]$ is an AF-domain. In
the present paper, we give partial results settling in the
affirmative the above question (Q). First, we prove that an
AF-domain $A$ is in fact a GAF-domain, thus in particular, any
finitely generated algebra over $k$ which is a domain is a
GAF-domain. Also, through Proposition 2.5, we prove that (Q) has a
positive answer in the case where $A$ is one-dimensional. Whereas,
our main result, Theorem 2.8, tackles the case $n=1$ of $(Q)$. It
computes dim$(A\otimes_kB)$ for a $k$-algebra $A$ such that $A[X]$
is an AF-domain and for an arbitrary $k$-algebra $B$ generalizing
Wadsworth's main theorem [16, Theorem 3.7] and further asserts that
$A$ is a GAF-domain. We end this paper by an example of a GAF-domain
$A$ such that, for any positive integer $n$, the polynomial ring
$A[n]$ is not an AF-domain.

Recent developments on height and grade of (prime) ideals as well as
on dimension theory in tensor products of $k$-algebras are to be
found in [2-6].\newpage

\noindent {\bf 2. Main results}\\

In this section, we handle the question (Q) set above. \\

First, for the convenience of the reader, we catalog some basic facts and results
connected with
the tensor product of
$k$-algebras. These will be used frequently in the sequel without explicit mention.\\

Let $A$ and $B$ be two $k$-algebras. If $p$ is a prime ideal of $A$,
$r=$ t.d.($\displaystyle {\frac Ap})$ and $\overline
{x_1},...,\overline {x_r}$ are elements of $\displaystyle {\frac
Ap}$, algebraically independent over $k$, with the $x_i \in A$, then
it is easily seen that $x_1,...,x_r$ are algebraically independent
over $k$ and $p\cap S=\emptyset$, where
$S=k[x_1,...,x_r]\setminus\{0\}$. If $A$ is an integral domain, then
$ht(p)+$t.d.($\displaystyle {\frac Ap})\leq$ t.d.($A)$ for each
prime ideal $p$ of $A$ (cf. [14, p. 37] ). Now, assume that $S_1$
and $S_2$ are multiplicative subsets of $A$ and $B$, respectively,
then $S_1^{-1}A\otimes _kS_2^{-1}B\cong S^{-1}(A\otimes _kB)$, where
$S=\{s_1\otimes s_2:s_1\in S_1$ and $s_2\in S_2\}$. We assume
familiarity with the natural isomorphisms for tensor products. In
particular, we identify $A$ and $B$ with their respective images in
$A\otimes_kB$. Also, $A\otimes_kB$ is a free (hence faithfully flat)
extension of $A$ and $B$. Moreover, recall that an AF-domain $A$ is
a locally Jaffard domain, that is, $ht(p[n])=ht(p)$ for each prime
ideal $p$ and each positive integer $n$ [16, Corollary 3.2].
Finally, we refer the reader to the useful result of Wadsworth [16,
Proposition 2.3] which yields a classification of the prime
ideals of $A\otimes_kB$ according to their contractions to $A$ and $B$.\\

We begin by recalling from [2], [5] and [16] the following useful results.\\

Our first result allows to construct a bunch of $k$-algebras $A$
arising from pullbacks which are not AF-domains while there exists
an integer $n\geq 1$ such that $A[n]$ is an AF-domain.\\

\noindent {\bf Proposition 2.1 [5, Proposition 2.2].} {\it Let $T$
be an integral domain which is a $k$-algebra, $M$ a maximal ideal of
$T$, $K:=\displaystyle {\frac TM}$ and $\varphi: T\rightarrow K$ the
canonical surjective homomorphism. Let $D$ be a proper subring of
$K$ and $A:=\varphi^{-1}(D)$. Assume that $T$ and $D$ are
AF-domains. Let $r:=$ t.d.$(K:k)$ and $s=$ t.d.$(D:k)$. Then the
polynomial ring
$A[r-s]$ is an AF-domain.}\\

Recall that, by [9], under the hypotheses of Proposition 2.1, $A$ is
an AF-domain if and only if t.d.($K:D)=r-s=0$. Thus, whenever $r>s$,
the issued pullback $A$ is
not an AF-domain.\\

\noindent {\bf Proposition 2.2 [2, Lemma 1.3].} {\it Let $A$ and $B$
be $k$-algebras such that $B$ is a domain. Let $p$ be a prime ideal
of $A$. Then, for each prime ideal $P$ of $A\otimes_kB$ which is
minimal over $p\otimes_kB$,
$$ht(P)=ht(p\otimes_kB)=ht(p[\mbox {t.d.}(B)]).$$}

\noindent {\bf Proposition 2.3 [16, Proposition 2.3].} {\it Let $A$
and $B$ be $k$-algebras and let $p\subseteq p^{\prime}$ be prime
ideals of $A$ and $q\subseteq q^{\prime}$ be prime ideals of $B$.
Then the natural ring homomorphism $\varphi: \displaystyle {\frac
{A\otimes_kB}{p\otimes_kB+A\otimes_kq}\longrightarrow \frac
Ap\otimes_k\frac Bq}$ such that $\varphi (\overline
{a\otimes_kb})=\overline a\otimes_k\overline b$ for each $a\in A$
and each $b\in B$, is an isomorphism and $$\varphi (\displaystyle
{\frac
{p^{\prime}\otimes_kB+A\otimes_kq^{\prime}}{p\otimes_kB+A\otimes_kq})=\frac
{p^{\prime}}p\otimes_k\frac Bq+\frac Ap\otimes_k\frac {q^{\prime}}q}.$$}\\

We establish the following easy result which is probably well known
but we have not located
references in the literature.\\

\noindent {\bf Proposition 2.4.} {\it Let $A$ be ring. Let
$I\subseteq J$ be ideals in $A$. Then $$ht(I)+ht(\displaystyle
{\frac JI})\leq ht(J).$$}

\noindent {\bf Proof.} If both $I$ and $J$ are prime ideals, then
the result easily follows. Fix a prime ideal $Q$ of $A$ that
contains $J$. Let $P$ be a minimal prime ideal of $I$ contained in
$Q$. As $ht(I)\leq ht(P)$ and $ht(P)+ht(\displaystyle {\frac
QP})\leq ht(Q)$, we get $ht(I)+ht(\displaystyle {\frac QP})\leq
ht(Q)$ for each minimal prime ideal $P$ of $I$ contained in $Q$.
Hence $ht(I)+$max$\{ht(\displaystyle {\frac QP}):P$ is a minimal
prime ideal of $A$ over $I$ contained in $Q\}=ht(I)+ht(\displaystyle
{\frac QI})\leq ht(Q)$. Therefore $$ht(I)+ht(\displaystyle {\frac
JI})\leq ht(I)+ht(\displaystyle {\frac QI})\leq ht(Q)$$ for each
prime ideal $Q$ of $A$ containing $J$. It follows that\\
$ht(I)+ht(\displaystyle {\frac JI})\leq$ min$\{ht(Q):J\subseteq
Q\in$ Spec$(A)\}=ht(J)$, as desired.
$\Box$\\

We begin by proving that an AF-domain is a GAF-domain.\\

\noindent {\bf Proposition 2.5.} {\it Let $A$ be an AF-domain. Then
$A$ is a GAF-domain.}\\

\noindent {\bf Proof.} Let $P$ be a prime ideal of $A\otimes_kB$, $p=P\cap A$ and $q=P\cap B$.
 By the special chain theorem for tensor products [5, Lemma 1.5],
 $$ht(P)=ht(A\otimes_kq)+ht(\frac P{A\otimes_kq}).$$ Also, note that
$\displaystyle {\frac P{A\otimes_kq}}$ is a prime ideal of
$\displaystyle {A\otimes_k\frac Bq}$ such that $\displaystyle {\frac
P{A\otimes_kq}\cap \frac Bq}=(\overline 0)$. Then $\displaystyle
{\frac P{A\otimes_kq}}$ survives in the localization
$A\otimes_kk_B(q)$ of $A\otimes_k\displaystyle {\frac Bq}$, where
$k_B(q)$ denotes the quotient field of $\displaystyle {\frac Bq}$,
so that, by a second application of [5, Lemma 1.5], we get

$$
\begin{array}{lll}
ht(\displaystyle {\frac P{A\otimes_kq})}&=&\displaystyle
{ht(p\otimes_k\frac Bq)+ht\Big {(}\frac
{P/(A\otimes_kq)}{p\otimes_k(B/q)}\Big {)}}\\
&&\\
&=&ht(\displaystyle {\frac {p\otimes_kB+A\otimes_kq}{A\otimes_kq})+
ht(\frac P{p\otimes_kB+A\otimes_kq})},
\end{array}$$

\noindent as $\displaystyle {\frac
{p\otimes_kB+A\otimes_kq}{A\otimes_kq}}\cong \displaystyle
{p\otimes_k\frac Bq}$ via Proposition 2.3. Applying Proposition 2.4,
it follows that

$$\begin{array}{lll}
ht(P)&=&ht(A\otimes_kq)+ht(\displaystyle {\frac
{p\otimes_kB+A\otimes_kq}{A\otimes_kq})+
ht(\frac P{p\otimes_kB+A\otimes_kq})}\\
&&\\
&\leq&ht(p\otimes_kB+A\otimes_kq)+ht(\displaystyle {\frac P{p\otimes_kB+A\otimes_kq})}\\
&\leq& ht(P).
\end{array}
$$

\noindent Then the equality holds, completing the proof. $\Box$\\

The following result settles in the affirmative the above question
(Q) in the case where $A$ is one-dimensional.\\

\noindent {\bf Proposition 2.6.} {\it Let $A$ be a one-dimensional
domain such that $A[n]$ is an AF-domain for some positive integer
$n$. Then
$A$ is a GAF-domain.}\\

\noindent {\bf Proof.} Let $B$ be a $k-$algebra and $P$ a prime
ideal of $A\otimes_kB$ with $p=P\cap A$ and $q=P\cap B$. Then,
applying [2, Theorem 1.1], we get $$ht(P)=\mbox {max}\Big
{\{}ht(q_1[\mbox {t.d.}(A)])+ht(\displaystyle {\frac q{q_1}[\mbox
{t.d.}(\frac Ap)])}+ht(\displaystyle {p[\mbox {t.d.}(\frac
B{q_1})]):}$$ $$q_1\in \mbox {Spec}(B) \mbox { with } q_1\subseteq
q\Big {\}}+\displaystyle { ht(\frac P{p\otimes_kB+A\otimes_kq})}.$$
Then, for each minimal prime ideal $Q$ of $p\otimes_kB+A\otimes_kq$,
we get $$ht(Q)=\mbox {max}\Big {\{}ht(q_1[\mbox
{t.d.}(A)])+ht(\displaystyle {\frac q{q_1}[\mbox {t.d.}(\frac
Ap)])}+ht(\displaystyle {p[\mbox {t.d.}(\frac B{q_1})])}:$$
$$q_1\in \mbox {Spec}(B) \mbox { with } q_1\subseteq q\Big {\}},$$
so that\newpage
$$ht(p\otimes_kB+A\otimes_kq)=\mbox {max}\Big
{\{}ht(q_1[\mbox {t.d.}(A)])+ht(\displaystyle {\frac q{q_1}[\mbox
{t.d.}(\frac Ap)])}+ht(\displaystyle {p[\mbox {t.d.}(\frac
B{q_1})]):}$$ $$q_1\in \mbox {Spec}(B) \mbox { with } q_1\subseteq
q\Big {\}}.$$ It follows that
$$ht(P)=ht(p\otimes_kB+A\otimes_kq)+ht(\displaystyle {\frac
P{p\otimes_kB+A\otimes_kq}}).$$ Then $A$ is a GAF-domain, as
desired. $\Box$\\

Next, we announce the principal result of this paper. It tackles the
case $n=1$ of the above-sited question $(Q)$ and generalizes the
main theorem of Wadsworth in [16], namely: If $A$ and $B$ are
$k$-algebras such that $A$ is an AF-domain, then\\
$$\begin{array}{lll}
\mbox {dim}(A\otimes_kB)&=&D\Big {(}\mbox { t.d.}(A),\mbox { dim}(A),B\Big {)}\\
&:=&\mbox {max}\Big {\{}ht(q[\mbox {t.d.}(A)])+\mbox {min}\Big
{(}\mbox {t.d.}(A),\mbox { dim}(A)+\mbox {t.d.}(\displaystyle {\frac
Bq})\Big {)}:q\in \mbox {Spec}(B)\Big {\}}\\
&& \mbox { [16, Theorem 3.7]}. \end{array} $$ This equality might be
rewritten in the following way which evokes our next general
formula,
$$\begin{array}{lll} \mbox {dim}(A\otimes_kB)&=&\mbox {max}\Big
{\{}ht(q[\mbox {t.d.}(A)])+\displaystyle { ht(p[\mbox {t.d.}(\frac
Bq)])+\mbox {min}\Big {(}\mbox {t.d.}(\frac Ap)},\mbox
{t.d.}(\displaystyle
{\frac Bq})\Big {)}:\\
&&p\in \mbox { Spec}(A)\mbox { and }q\in \mbox {Spec}(B)\Big
{\}}(\mbox { as }A \mbox { is a locally Jaffard domain}).
\end{array}$$ First, it is worthwhile recalling the following
definition and results from [3] and [5]. Let $A$ and $B$ be
$k$-algebras and $P$ be a prime ideal of $A\otimes _kB$. Let
 $q_0\in$ Spec$(B)$ such that $q_0\subset P\cap B$. We denote by
$\lambda \Big {(}(.,q_0),P\Big {)}$ the maximum of lengths of chains
of prime ideals of $A\otimes_kB$ of the form $P_0\subset P_1\subset
...\subset P_s=P$ such that $P_i\cap B=q_0$, for $i=0,1,...,s-1$.
Applying [3, Lemma 2.4], if $A$ and $B$ are integral domains, then
$$\lambda {\Big (}(.,(0)),P{\Big )}\leq\mbox { t.d.}(A)-\mbox
{t.d.}(\frac Ap)+ht(q[\mbox {t.d.}(\frac Ap)])+ht(\frac
P{p\otimes_kB+A\otimes_kq}).$$ Further, recall that, if $A$ is a
$k$-algebra and $n\geq 0$ is an integer, then the polynomial ring
$A[n]$ is an AF-domain if and only if
$$ht(p[n])+\mbox {t.d.}\displaystyle {(\frac Ap)}=\mbox {
t.d.}(A)$$ for each prime
ideal $p$ of $A$ [5, Lemma 2.1].\\

\noindent {\bf Theorem 2.7.} {\it Let $A$ be a $k$-algebra such that
the polynomial ring $A[X]$ is an AF-domain. Let $B$ be an arbitrary
$k$-algebra. Then the following statements hold:\\

a) If $P$ is a prime ideal of $A\otimes_kB$, $p=P\cap A$ and
$q=P\cap B$, then\newpage $$ht(P)=\mbox { max}\Big {\{}ht(q_1[\mbox
{t.d.(A)}])+ht(p_1[\mbox {t.d.}(\displaystyle {\frac
B{q_1})])+ht(\frac q{q_1}[\mbox {t.d.}(\frac A{p_1})])+ht(\frac
p{p_1}[\mbox {t.d.}(\frac Bq)])}:$$ $$ \displaystyle { p_1\subseteq
p \mbox { and }q_1\subseteq q\mbox { are prime ideals of } A\mbox {
and }B,\mbox { respectively}}\Big {\}}+$$ \hspace
{2cm}$\displaystyle {ht(\frac P{p\otimes_kB+A\otimes_kq})}$.

b) dim$(A\otimes_kB)=\mbox { max}\Big {\{}ht(q_1[\mbox
{t.d.(A)}])+ht(p_1[\mbox {t.d.}(\displaystyle {\frac
B{q_1})])+ht(\frac q{q_1}[\mbox {t.d.}(\frac A{p_1})])}+$

$\displaystyle {ht(\frac p{p_1}[\mbox {t.d.}(\frac
Bq)])}+\displaystyle {\mbox {min\Big {(}t.d.(}\frac Ap),\mbox
{t.d.}(\frac Bq)\Big {)}}:p_1\subseteq p\in$ Spec$(A)$ and
$q_1\subseteq q\in$ Spec$(B)\Big {\}}.$\\

c) $A$ is a GAF-domain.}\\

\noindent {\bf Proof.} a) {\bf  \underline {Step 1}.} $B$ is an integral domain.\\

\noindent If t.d.$(B)=0$, then $B$ is an algebraic extension field
of $k$, and thus, by [5, Lemma 1.5], we are done. Next, assume that
t.d.$(B)\geq 1$. By Proposition 2.2 and Proposition 2.3, we get,
$\forall p_1\subseteq p\in$ Spec$(A)$ and $\forall q_1\subseteq
q\in$
Spec$(B)$,\\

$\begin{array}{lll} ht(q_1[\mbox
{t.d.}(A)])&=&ht(A\otimes_kq_1),\\
&&\\
\displaystyle {ht(p_1[\mbox {t.d.}(\frac B{q_1})])}&=&\displaystyle {ht(p_1\otimes_k\frac B{q_1})}=
\displaystyle {ht(\frac {p_1\otimes_kB+A\otimes_kq_1}{A\otimes_kq_1})},\\
&&\\
\displaystyle {ht(\frac q{q_1}[\mbox {t.d.}(\frac
A{p_1})])}&=&\displaystyle {ht(\frac
{p_1\otimes_kB+A\otimes_kq}{p_1\otimes_kB+A\otimes_kq_1})}\mbox {, and}\\
&&\\
\displaystyle {ht(\frac p{p_1}[\mbox {t.d.}(\frac
B{q})])}&=&\displaystyle {ht(\frac
{p\otimes_kB+A\otimes_kq}{p_1\otimes_kB+A\otimes_kq}).}
\end{array}$\\

\noindent It follows that, $\forall p_1\subseteq p\in$ Spec$(A)$
and $\forall q_1\subseteq q\in$ Spec$(B)$,\\

$ht(q_1[\mbox {t.d.(A)}])+ht(p_1[\mbox {t.d.}(\displaystyle {\frac
B{q_1})])+ht(\frac q{q_1}[\mbox {t.d.}(\frac A{p_1})])+ht(\frac
p{p_1}[\mbox {t.d.}(\frac Bq)])}=$\\

$ht(A\otimes_kq_1)+\displaystyle {ht(\frac
{p_1\otimes_kB+A\otimes_kq_1}{A\otimes_kq_1})}+\displaystyle
{ht(\frac
{p_1\otimes_kB+A\otimes_kq}{p_1\otimes_kB+A\otimes_kq_1})}+\displaystyle
{ht(\frac
{p\otimes_kB+A\otimes_kq}{p_1\otimes_kB+A\otimes_kq})\leq}$
$$ht(p\otimes_kB+A\otimes_kq) \mbox {, by Proposition 2.4}.$$\\
Consequently, $$\mbox { max}\Big {\{}ht(q_1[\mbox
{t.d.(A)}])+ht(p_1[\mbox {t.d.}(\displaystyle {\frac
B{q_1})])+ht(\frac q{q_1}[\mbox {t.d.}(\frac A{p_1})])+ht(\frac
p{p_1}[\mbox {t.d.}(\frac Bq)])}:\displaystyle { p_1\subseteq p
\mbox { and }}$$ $$q_1\subseteq q \mbox { are prime ideals of }
A\mbox { and }B,\mbox { respectively}\Big {\}}+ht(\displaystyle
{\frac P{p\otimes_kB+A\otimes_kq})}\leq$$
$$ht(p\otimes_kB+A\otimes_kq)+ht(\displaystyle {\frac
P{p\otimes_kB+A\otimes_kq})}\leq ht(P).$$ Our proof of the reverse
inequality uses induction on $ht(p)$ and $ht(q)$. First, note that\\

(1)\hspace {1cm} $\mbox {max}\Big {\{}ht(q[\mbox
{t.d.}(A)])+ht(p[\displaystyle {\mbox {t.d.}(\frac Bq)]),ht(p[\mbox
{t.d.}(B)])+ ht(q[\mbox {t.d.}(\frac Ap)])\Big {\}}}\leq $\\ $$\mbox
{ max}\Big {\{}ht(q_1[\mbox {t.d.(A)}])+ht(p_1[\mbox
{t.d.}(\displaystyle {\frac B{q_1})])+ht(\frac q{q_1}[\mbox
{t.d.}(\frac A{p_1})])+ht(\frac p{p_1}[\mbox {t.d.}(\frac
Bq)])}:\displaystyle { p_1\subseteq p \mbox { and }}$$
$$q_1\subseteq q \mbox { are prime ideals of } A\mbox { and }B,\mbox
{ respectively}\Big {\}},\mbox {  it suffices to take}$$ $$\left \{\begin{array}{lll}p_1=p\\
\\q_1=q  \end{array}\right.\mbox { and }\left
\{\begin{array}{lll}p_1=p\\
\\q_1=(0). \end{array}\right.$$ The case where either $ht(p)=0$ or
$ht(q)=0$ is fairly easy applying [5, Lemma 1.5]. Then, assume that
$ht(p)>0$ and $ht(q)>0$. Suppose that t.d.$(\displaystyle {\frac
Bq})\geq 1$ and let $x\in B$ such that $\overline x$ is a
transcendental element of $\displaystyle {\frac Bq}$ over $k$, and
put $S:=k[x]\setminus \{0\}$. Then,\\

$\left \{ \begin{array}{ll} x$ is transcendental over $A$ (we
identify $x$ with its image $1\otimes_kx$ through the$\\
\hspace{4.5cm} $canonical injection $B\rightarrow A\otimes_kB)&\\
&\\
S^{-1}(A\otimes_kB)\cong A\otimes_kS^{-1}B\cong \Big
{(}A\otimes_kk(x)\Big {)}\otimes_{k(x)}S^{-1}B\cong
S^{-1}A[x]\otimes_{k(x)}S^{-1}B\\
&\\
q\cap S=\emptyset\\
&\\
S^{-1}P\in$ Spec$\Big {(}S^{-1}A[x]\otimes_{k(x)}S^{-1}B\Big {)}\\
&\\
S^{-1}A[x]$ is an AF-domain, by hypotheses.$
\end{array}
\right.$\\

\noindent Hence
$$\begin{array}{lll}
ht(P)&=&ht(S^{-1}P)=ht\Big {(}S^{-1}A[x]\otimes_{k(x)}S^{-1}q\Big
{)}+\displaystyle
{ht\Big {(}\frac {S^{-1}P}{S^{-1}A[x]\otimes_{k(x)}S^{-1}q}\Big {)}},\\
&&\mbox { via [5, Lemma 1.5]}\\
&=&ht\Big {(}(A\otimes_kk(x))\otimes_{k(x)}S^{-1}q\Big
{)}+\displaystyle {ht\Big {(}\frac
{S^{-1}P}{(A\otimes_kk(x))\otimes_{k(x)}S^{-1}q}\Big {)}}\\
&=&ht(A\otimes_kS^{-1}q)+\displaystyle {ht(\frac
{S^{-1}P}{A\otimes_kS^{-1}q})}
\end{array}$$
$$\begin{array}{lll}
&=&ht(A\otimes_kq)+\displaystyle {ht(\frac P{A\otimes_kq})}\\
&=&\displaystyle { ht(q[\mbox {t.d.}(A)])+ht(p[\mbox { t.d.}(\frac
Bq)])+ht(\frac P{p\otimes_kB+A\otimes_kq})}, \mbox { by [5, Lemma
1.5]},
\end{array}
$$
$$\begin{array}{lll}
&&\mbox { since }\displaystyle {\frac P{A\otimes_kq}\cap \frac
Bq=(0)}\mbox { and }\displaystyle {p\otimes_k\frac Bq}\cong
\displaystyle {\frac {p\otimes_kB+A\otimes_kq}{A\otimes_kq}}\mbox {
via Proposition 2.3,}\\
&&\\
&&\mbox { so that } \displaystyle {\frac
{P/(A\otimes_kq)}{p\otimes_k(B/q)}\cong
\frac P{p\otimes_kB+A\otimes_kq}}\\
\end{array}$$
$$\begin{array}{lll}
&\leq& \mbox { max}\Big {\{}ht(q_1[\mbox {t.d.(A)}])+ht(p_1[\mbox
{t.d.}(\displaystyle {\frac B{q_1})])+ht(\frac q{q_1}[\mbox
{t.d.}(\frac A{p_1})])+ht(\frac p{p_1}[\mbox {t.d.}(\frac
Bq)])}:\\
&&\displaystyle { p_1\subseteq p \mbox { and }}q_1\subseteq q \mbox
{ are prime ideals of } A\mbox { and }B,\mbox { respectively}\Big
{\}}+\\
&&\displaystyle {ht(\frac P{p\otimes_kB+A\otimes_kq})},\mbox { by } (1)\\
&&\\
&\leq&ht(P),\mbox { then the equality holds, as contended}.
\end{array}$$
Next, suppose that t.d.$(\displaystyle {\frac Bq)=0}$. Then, by [16,
Proposition 2.3], $P$ is a minimal prime ideal of
$p\otimes_kB+A\otimes_kq$. Let $Q\in$ Spec$(A\otimes_kB)$ such that
$Q\subset P$ and $ht(P)=1+ht(Q)$.  Let $p^{\prime}=Q\cap A$ and
$q^{\prime}=Q\cap B$. Then either $p^{\prime}\subset p$ or
$q^{\prime}\subset q$. Three cases arise.\\

\noindent {\bf Case 1. $q^{\prime}\subset q$ and $q^{\prime}\neq
(0)$.} Then t.d.$(\displaystyle {\frac B{q^{\prime}})}\geq 1$ since
$1\leq ht(\displaystyle {\frac q{q^{\prime}})}+$t.d.$(\displaystyle
{\frac {B/q^{\prime}}{q/q^{\prime}})}\leq$ t.d.$(\displaystyle
{\frac B{q^{\prime}})}$, so that, by the above discussion,
$$ht(Q)=ht(A\otimes_kq^{\prime})+ht(\displaystyle {\frac
Q{A\otimes_kq^{\prime}})},\mbox { and hence }
ht(P)=ht(A\otimes_kq^{\prime})+ht(\displaystyle {\frac
P{A\otimes_kq^{\prime}})}.$$ As $q^{\prime}\neq (0)$,
$ht(\displaystyle {\frac q{q^{\prime}})<ht(q)}$, then by inductive
hypotheses with respect to $A\otimes_k\displaystyle {\frac
B{q^{\prime}}}$, we get
$$ht(\displaystyle {\frac P{A\otimes_kq^{\prime}})=\mbox { max}\Big {\{}ht(\frac {q_1}{q^{\prime}}
[\mbox {t.d.(A)}])+
ht(p_1[\mbox {t.d.}}(\displaystyle {\frac B{q_1})])+ht(\frac
q{q_1}[\mbox {t.d.}(\frac A{p_1})])+ht(\frac p{p_1}[\mbox
{t.d.}(\frac Bq)])}:$$
$$ \displaystyle { q^{\prime}\subseteq q_1\subseteq q \mbox { and }p_1\subseteq p\mbox {
are prime ideals of } A\mbox { and }B,\mbox { respectively}}\Big
{\}}+$$ \hspace {2cm}$\displaystyle {ht(\frac
P{p\otimes_kB+A\otimes_kq})}$, as, by Proposition 2.3,
$$\displaystyle {\frac
{P/(A\otimes_kq^{\prime})}{p\otimes_k(B/q^{\prime})+A\otimes_k(q/q^{\prime})}\cong
\frac
{P/(A\otimes_kq^{\prime})}{(p\otimes_kB+A\otimes_kq)/(A\otimes_kq^{\prime})}\cong
\frac P{p\otimes_kB+A\otimes_kq}}.$$ Therefore
$$
\begin{array}{lll}
ht(P)&=&\displaystyle {ht(q^{\prime}[\mbox {t.d.}(A)])+\mbox
{max}\Big {\{}ht(\frac {q_1}{q^{\prime}}[\mbox {t.d.(A)}])+
ht(p_1[\mbox
{t.d.}(\frac B{q_1})])+}\\
&&\displaystyle {ht(\frac q{q_1}[\mbox {t.d.}(\frac
A{p_1})])+ht(\frac p{p_1}[\mbox {t.d.}(\frac Bq)])}:p_1\subseteq
p\mbox { and
}q^{\prime}\subseteq q_1\subseteq q \mbox { are prime ideals of }\\
&&\displaystyle { A\mbox { and }B,\mbox { respectively}}\Big
{\}}+\displaystyle {ht(\frac P{p\otimes_kB+A\otimes_kq})}
\end{array}$$
$$\begin{array}{lll}
 &=&\displaystyle {\mbox { max}\Big {\{}ht(q^{\prime}[\mbox
{t.d.}(A)])+ht(\frac {q_1}{q^{\prime}}[\mbox {t.d.(A)}])+
ht(p_1[\mbox
{t.d.}(\frac B{q_1})])+}\\
&&\displaystyle {ht(\frac q{q_1}[\mbox {t.d.}(\frac
A{p_1})])+ht(\frac p{p_1}[\mbox {t.d.}(\frac Bq)])}:p_1\subseteq
p\mbox { and
}q^{\prime}\subseteq q_1\subseteq q \mbox { are prime ideals of }\\
&&\displaystyle { A\mbox { and }B,\mbox { respectively}}\Big
{\}}+\displaystyle {ht(\frac
P{p\otimes_kB+A\otimes_kq})}\\
&&\\
&\leq&\displaystyle {\mbox { max}\Big {\{}ht(q_1[\mbox {t.d.}(A)])+
ht(p_1[\mbox {t.d.}(\frac B{q_1})])+\displaystyle {ht(\frac
q{q_1}[\mbox {t.d.}(\frac
A{p_1})])+ht(\frac p{p_1}[\mbox {t.d.}(\frac Bq)])}:}\\
&&p_1\subseteq p\mbox { and }q_1\subseteq q\mbox { are prime ideals
of }\displaystyle { A\mbox { and }B,\mbox {
respectively}}\Big {\}}+\\
&& \displaystyle {ht(\frac
P{p\otimes_kB+A\otimes_kq})}\\
&\leq&ht(P).

\end{array}$$ Then the equality holds.\\

\noindent {\bf Case 2. $q^{\prime}=q$.} Then $p^{\prime}\subset p$.
By inductive hypotheses, we get
$$ht(Q)=\mbox { max}\Big {\{}ht(q_1[\mbox {t.d.(A)}])+ht(p_1[\mbox
{t.d.}(\displaystyle {\frac B{q_1})])+ht(\frac q{q_1}[\mbox
{t.d.}(\frac A{p_1})])+ht(\frac {p^{\prime}}{p_1}[\mbox {t.d.}(\frac
Bq)])}:$$
$$ \displaystyle { p_1\subseteq p^{\prime}\mbox { and }q_1\subseteq q \mbox {
are prime ideals of } A\mbox { and }B,\mbox { respectively}}\Big
{\}}+$$ \hspace {2cm}$\displaystyle {ht(\frac
Q{p^{\prime}\otimes_kB+A\otimes_kq})}.$ Hence $$ht(P)\leq\mbox {
max}\Big {\{}ht(q_1[\mbox {t.d.(A)}])+ht(p_1[\mbox
{t.d.}(\displaystyle {\frac B{q_1})])+ht(\frac q{q_1}[\mbox
{t.d.}(\frac A{p_1})])+ht(\frac {p^{\prime}}{p_1}[\mbox {t.d.}(\frac
Bq)])}:$$
$$ \displaystyle { p_1\subseteq p^{\prime}\mbox { and }q_1\subseteq q \mbox {
are prime ideals of } A\mbox { and }B,\mbox { respectively}}\Big
{\}}+$$ \hspace {2cm}$\displaystyle {ht(\frac
P{p^{\prime}\otimes_kB+A\otimes_kq})}.$\\ As $\displaystyle {\frac
P{p^{\prime}\otimes_kB+A\otimes_kq}\cap \frac Bq=(\overline 0)}$ and
$\displaystyle {\frac P{p^{\prime}\otimes_kB+A\otimes_kq}\cap \frac
A{p^{\prime}}=\frac p{p^{\prime}}}$, we get, by [5, Lemma 1.5],
$$\begin{array}{lll}
\displaystyle {ht(\frac
P{p^{\prime}\otimes_kB+A\otimes_kq})}&=&\displaystyle {ht(\frac
p{p^{\prime}}[\mbox {t.d.}(\frac Bq)])+ht\Big {(}\frac
{P/(p^{\prime}\otimes_kB+A\otimes_kq)}{
(p/p^{\prime})\otimes_k(B/q)}}\Big {)}\\
&&\\
&=&\displaystyle {ht(\frac p{p^{\prime}}[\mbox {t.d.}(\frac
Bq)])+ht(\frac P{p\otimes_kB+A\otimes_kq})} \end{array}$$ since
$\displaystyle {\frac p{p^{\prime}}\otimes_k\frac Bq}\cong
\displaystyle {\frac
{p\otimes_kB+A\otimes_kq}{p^{\prime}\otimes_kB+A\otimes_kq}}$ via
Proposition 2.3. It follows that
$$\begin{array}{lll}
ht(P)&\leq&\mbox { max}\Big {\{}ht(q_1[\mbox
{t.d.(A)}])+ht(p_1[\mbox {t.d.}(\displaystyle {\frac
B{q_1})])+ht(\frac q{q_1}[\mbox {t.d.}(\frac A{p_1})])+ht(\frac
{p^{\prime}}{p_1}[\mbox {t.d.}(\frac
Bq)])}:\\
&& \displaystyle {p_1\subseteq p^{\prime}\mbox { and } q_1\subseteq
q\mbox { are prime ideals of } A\mbox { and }B,\mbox {
respectively}}\Big {\}}+\\
&&\displaystyle {ht(\frac p{p^{\prime}}[\mbox {t.d.}(\frac
Bq)])+ht(\frac P{p\otimes_kB+A\otimes_kq})}\\
&&\\
&\leq&\mbox { max}\Big {\{}ht(q_1[\mbox {t.d.(A)}])+ht(p_1[\mbox
{t.d.}(\displaystyle {\frac B{q_1})])+ht(\frac q{q_1}[\mbox
{t.d.}(\frac A{p_1})])+ht(\frac {p^{\prime}}{p_1}[\mbox {t.d.}(\frac
Bq)])}+\\
&& \displaystyle {ht(\frac p{p^{\prime}}[\mbox {t.d.}(\frac
Bq)])}:\displaystyle { p_1\subseteq
p^{\prime}\mbox { and }q_1\subseteq q\mbox { are prime ideals of } A\mbox { and }B,}\\
&&\mbox {
respectively}\Big {\}}+\displaystyle {ht(\frac P{p\otimes_kB+A\otimes_kq})}\\
&&\\
&\leq&\mbox { max}\Big {\{}ht(q_1[\mbox {t.d.(A)}])+ht(p_1[\mbox
{t.d.}(\displaystyle {\frac B{q_1})])+ht(\frac q{q_1}[\mbox
{t.d.}(\frac A{p_1})])+ht(\frac p{p_1}[\mbox {t.d.}(\frac
Bq)])}:\\
&& p_1\subseteq p\mbox { and }q_1\subseteq q\mbox { are prime ideals
of } A\mbox { and }B,\mbox {
respectively}\Big {\}}+\\
&&\displaystyle {ht(\frac P{p\otimes_kB+A\otimes_kq})}\\
&\leq& ht(P), \mbox { and then the equality holds.}

\end{array}$$

\noindent {\bf Case 3. $q^{\prime}=(0)$.} Then\\

$\begin{array}{lll} ht(P)&=&1+ht(Q)=\lambda \Big {(}(.,(0)),P\Big {)}\\
&\leq&\mbox {t.d.}(A)-\displaystyle {\mbox {t.d.}(\frac
Ap)+ht(q[\mbox {t.d.}(\frac
Ap)])+ht(\frac P{p\otimes_kB+A\otimes_kq})},\mbox { by [3, Lemma 2.4]}\\
&=&\displaystyle {ht(p[X])+ht(q[\mbox {t.d.}(\frac Ap)])+ht(\frac
P{p\otimes_kB+A\otimes_kq})},$ as $A[X]$ is an AF-domain$\\
&\leq& \displaystyle {ht(p[\mbox {t.d.}(B)])+ht(q[\mbox {t.d.}(\frac
Ap)])+ht(\frac P{p\otimes_kB+A\otimes_kq})},\mbox { as t.d.}(B)\geq
1
\end{array}$
$\begin{array}{lll}
&\leq&\mbox { max}\Big
{\{}ht(q_1[\mbox {t.d.(A)}])+ht(p_1[\mbox {t.d.}(\displaystyle
{\frac B{q_1})])+ht(\frac q{q_1}[\mbox {t.d.}(\frac
A{p_1})])+ht(\frac p{p_1}[\mbox {t.d.}(\frac
Bq)])}:\\
&& p_1\subseteq p\in\mbox { Spec}(A)\mbox { and }q_1\subseteq q\in\mbox { Spec}(B)\Big {\}}+\\
&&\displaystyle {ht(\frac P{p\otimes_kB+A\otimes_kq})},\mbox { by }(1)\\
&\leq& ht(P).
\end{array}$\\ Then the equality holds, as desired.\\

\noindent {\bf \underline {Step 2}.} $B$ is an arbitrary $k$-algebra.\\

\noindent Let $P_0\subset P_1\subset...\subset P_h=P$ be a chain of
prime ideals of $A\otimes_kB$ such that $h=ht(P)$. Let $q_0:=P_0\cap
B$. Then
$$\displaystyle {\frac {P_0}{A\otimes_kq_0}\subset \frac
{P_0}{A\otimes_kq_0}\subset \frac {P_1}{A\otimes_kq_0}\subset
...\subset \frac {P_h}{A\otimes_kq_0}=\frac {P}{A\otimes_kq_0}}$$ is
a chain of prime ideals of $A\otimes_k\displaystyle {\frac B{q_0}}$
and $h=ht(P)=ht(\displaystyle {\frac P{A\otimes_kq_0})}$. By Step 1,
$$\begin{array}{lll}
ht(P)&=&ht(\displaystyle {\frac P{A\otimes_kq_0})=\mbox { max}\Big
{\{}ht(\frac {q_1}{q_0}}[\mbox {t.d.(A)}])+ht(p_1[\mbox
{t.d.}(\displaystyle {\frac B{q_1})])+ht(\frac q{q_1}[\mbox
{t.d.}(\frac A{p_1})])}+\\
&&\\
&&\displaystyle {ht(\frac p{p_1}[\mbox {t.d.}(\frac
Bq)])}:p_1\subseteq p\in\mbox { Spec}(A)\mbox { and
}q_0\subseteq q_1\subseteq q\in\mbox { Spec}(B)\Big {\}}+\\
&&\\
&&\displaystyle {ht(\frac {P/(A\otimes_kq_0)}{p\otimes_k
(B/q_0)+A\otimes_k (q/q_0)})}\\
&&\\
&\leq&\mbox { max}\Big {\{}ht(q_1[\mbox
{t.d.(A)}])+ht(p_1[\mbox {t.d.}(\displaystyle {\frac
B{q_1})])+ht(\frac q{q_1}[\mbox
{t.d.}(\frac A{p_1})])}+\\
&&\\
&&\displaystyle {ht(\frac p{p_1}[\mbox {t.d.}(\frac
Bq)])}:p_1\subseteq p\in\mbox { Spec}(A)\mbox { and
}q_1\subseteq q\in\mbox { Spec}(B)\Big {\}}+\\
&&\\
&&\displaystyle {ht(\frac {P}{p\otimes_k B+A\otimes_k q})}\\
&\leq&ht(P),\mbox { then the equality holds establishing the desired
formula}.
\end{array}$$

\noindent b) It is a direct consequence of (a) and [16, Proposition
2.3].\\

\noindent c) Let $B$ be a $k$-algebra. Let $p\in$ Spec$(A)$ and
$q\in$ Spec$(B)$. Applying (a), we get\newpage
$$ht(P)=\mbox { max}\Big {\{}ht(q_1[\mbox {t.d.(A)}])+ht(p_1[\mbox
{t.d.}(\displaystyle {\frac B{q_1})])+ht(\frac q{q_1}[\mbox
{t.d.}(\frac A{p_1})])+ht(\frac p{p_1}[\mbox {t.d.}(\frac Bq)])}:$$
$$ \displaystyle { p_1\subseteq p\mbox { and }q_1\subseteq q\mbox {
are prime ideals of } A\mbox { and }B,\mbox { respectively}}\Big
{\}}$$ for each minimal prime ideal $P$ of
$p\otimes_kB+A\otimes_kq$. It follows that
$$ht(p\otimes_kB+A\otimes_kq)=\mbox { max}\Big {\{}ht(q_1[\mbox
{t.d.(A)}])+ht(p_1[\mbox {t.d.}(\displaystyle {\frac
B{q_1})])+ht(\frac q{q_1}[\mbox {t.d.}(\frac A{p_1})])+}$$ $$
\displaystyle {ht(\frac p{p_1}[\mbox {t.d.}(\frac Bq)]):
p_1\subseteq p\mbox { and }q_1\subseteq q\mbox { are prime ideals of
} A\mbox { and }B,\mbox { respectively}}\Big {\}},$$ and thus
$$ht(P)=ht(p\otimes_B+A\otimes_kq)+ht(\displaystyle {\frac P{p\otimes_kB+A\otimes_kq}})$$ for each prime ideal
$P$ of $A\otimes_kB$ such that $p=P\cap A$ and $q=P\cap B$. Then $A$
is a GAF-domain, as desired. $\Box$\\

Next, we present an example of a GAF-domain $A$ such that $A[n]$
fails to be an AF-domain for any positive integer $n$.\\

\noindent {\bf Example 2.8.} Let $k$ be an algebraically closed
field. Consider the $k$-algebra homomorphism
$\varphi:k[X,Y]\rightarrow k[[t]]$ such that $\varphi (X)=t$ and
$\varphi (Y)=s:=\displaystyle {\sum_{n\geq 1}}t^{n!}$. Since $s$ is
known to be transcendental over $k(t)$, $\varphi$ is injective. This
induces an embedding $\overline {\varphi}:k(X,Y)\rightarrow k((t))$
of fields. Put $A=\overline {\varphi}^{-1}(k[[t]])$. It is easy to
check that $A$ is a discrete rank-one valuation overring of $k[X,Y]$
of the forme $A=k+p$, where $p=XA$. Note that, for each positive
integer $n$, and since $A$ is Noetherian, thus a locally Jaffard
domain,
$$ht(p[n])+\mbox {t.d.}(\displaystyle {\frac Ap})=ht(p)+\mbox
{t.d.}(\displaystyle {\frac Ap})=1<2=\mbox { t.d.}(A).$$ Then, via
[5, Lemma 2.1], for each positive integer $n$, $A[n]$ is not an
AF-domain. Let $B$ be an arbitrary $k$-algebra. Let $P$ be a prime
ideal of $A\otimes_kB$ and $q=P\cap B$. If $P\cap A=(0)$, then $P$
survives in $k(X,Y)\otimes_kB$, and thus, by [5, Lemma 1.5], we are
done. Now, assume that $P\cap A=p$. Then, as t.d.$(\displaystyle
{\frac Ap)=0}$, $P$ is minimal over $p\otimes_kB+A\otimes_kq$ [16,
Proposition 2.3]. Moreover, since $k$ is algebraically closed,
$p\otimes_kB+A\otimes_kq$ is a prime ideal of $A\otimes_kB$, as
$\displaystyle {\frac {A\otimes_kB}{p\otimes_kB+A\otimes_kq}\cong
\frac Ap\otimes_k\frac Bq}$ is an integral domain, by [17, Corollary
1, p. 198]. Hence $P=p\otimes_kB+A\otimes_kq$, so that
$ht(P)=ht(p\otimes_kB+A\otimes_kq)$. It follows that $A$ is a
GAF-domain, as desired. $\Box$\newpage

\noindent {\bf References}
\begin{list}{}{\topsep=2mm \parsep=0mm \itemsep=0mm \leftmargin=7.1mm}
{\small
\item [{[1]}]  J.T. Arnold, {\it On the dimension theory of
overrings of an integral domain}, Trans. Amer. Math. Soc., 138,
1969,313-326.
\item [{[2]}]  S. Bouchiba, {\it On Krull dimension of tensor products of algebras arising from AF-domains},
J. Pure Appl. Algebra 203 (2005) 237-251.
\item [{[3]}]  S. Bouchiba, {\it On chains of prime ideals of tensor products of algebras},
J. Pure Appl. Algebra 203 (2007) 237-251.
\item [{[4]}]  S. Bouchiba, {\it AF-domains and Jaffard's special
chain theorem}, The Arabian Journal for Science and Engineering,
Vol. 30, n: 2A, (2005), 1-7.
\item [{[5]}]  S. Bouchiba, F. Girolami, and S. Kabbaj, {\it The dimension of tensor
products of k-algebras arising from pullbacks}, J. Pure Appl.
Algebra 137 (1999), 125-138.
\item [{[6]}]  S. Bouchiba and S. Kabbaj, {\it Tensor products of Cohen-Macaulay rings: Solution to a problem of
Grothendieck}, J. Algebra 252 (2002) 65-73.
\item [{[7]}] J.W. Brewer, P.R. Montgomery, E.A. Rutter, W.J.
Heinzer, {\it Krull dimension of polynomial rings}, Lecture Notes in
Math., vol. 311, Springer-Verlag, Berlin, NY, 1972, 26-45.
\item [{[8]}]  R. Gilmer, {\it Multiplicative Ideal Theory}, Marcel Dekker, New
York, 1972.
\item [{[9]}]  F. Girolami, {\it AF-Rings and Locally Jaffard
Rings}, Lecture Notes Pure Appl. Math., Vol. 153. New York: Dekker,
1994, 151–161.
\item [{[10]}]  A. Grothendieck, {\it El\'ements de g\'eom\'etrie
alg\'ebrique}, Vol. IV,  Institut des Hautes Etudes Sci. Publ. Math.
No. 24, Bures-sur-yvette, 1965.
\item [{[11]}] P. Jaffard, {\it Th\'eorie de la dimension dans les anneaux de polyn\^omes}, M\'em. Sc. Math.
146, Gauthier-Villars, Paris, 1960.
\item [{[12]}] I. Kaplansky, {\it Commutative rings}, Chicago univ.,
Chicago press, 1974.
\item [{[13]}] H. Matsumura, {\it Commutative Ring Theory}, Cambridge University
Press, 1986.
\item [{[14]}] M. Nagata, {\it Local rings}, Interscience, New York,
1962.
\item [{[15]}]  R.Y. Sharp, {\it The dimension of the tensor product of two field
extensions}, Bull. London Math. Soc. 9 (1977), 42-48.
\item [{[16]}]  A.R. Wadsworth, {\it The Krull dimension of tensor products of
commutative algebras over a field}, J. London Math. Soc. 19 (1979),
391-401.
\item [{[17]}]  O. Zariski and P. Samuel, {\it Commutative Algebra Vol. I}, Van
Nostrand, Princeton, 1960.}

\end{list}

\end{document}